\documentclass[11pt]{article}
\usepackage[english]{babel}
\usepackage{amssymb}
\newtheorem{lemma}{Lemma}[section]

\newtheorem{theorem}[lemma]{Theorem}

\newtheorem{corollary}[lemma]{Corollary}
\newtheorem{example}[lemma]{Example}

\newtheorem{note}[lemma]{Remark}

\def\endproof{\hfill$\Box$}

\def\endproof{\hfill$\Box$}

\newcommand{\iso}{\,\mbox{\raisebox{0.65\height}{$\sim$}}\!\!\!\!\!%
\mbox{\raisebox{-0.4\height}{$\to$}}\,}

\title{Spherical orders, properties and countable spectra of
their theories\footnote{The research is supported by Committee of
Science in Education and Science Ministry of the Republic of
Kazakhstan, Grant No. AP08855544 (Sections \ref{sec3} and
\ref{sec4}) and Russian Scientific Foundation, Project
No.~22-21-00044 (Section \ref{sec6}). The work of the second
author was carried out in the framework of the State Contract of
the Sobolev Institute of Mathematics, Project No.~FWNF-2022-0012
(Section \ref{sec5}).}}
\author{B.Sh. Kulpeshov, S.V. Sudoplatov}

\date{}

\begin{document}

\maketitle
\begin{abstract}
We study semantic and syntactic properties of spherical orders and
their elementary theories, including finite and dense orders and
their theories. It is shown that theories of dense $n$-spherical
orders are countably categorical and decidable. The values for
spectra of countable models of unary expansions of $n$-spherical
theories are described. The Vaught conjecture is confirmed for
countable constant expansions of dense $n$-spherical theories.
\end{abstract}


\maketitle

\section{Introduction}\label{sec1}

Spherical orders \cite{aafot, ananat} are natural generalizations
of linear and circular orders \cite{KulMac, Kulp3, Kulp4} allowing
to compare tuples and their elements on $n$-spheres. These orders
present geometric realizations of arbitrary arities of formulae
along with known algebraic ones \cite{Semenov2003}.

In the present paper we study some kinds of $n$-spherical orders
including finite and dense ones, and their elementary theories. We
show that theories of dense $n$-spherical orders are countably
categorical and decidable. Basing on classical Ehrenfeucht
examples, Ehrenfeucht expansions of dense $n$-spherical orders are
constructed producing theories with arbitrarily finitely many
countable models. The Vaught conjecture is confirmed for countable
constant expansions of dense $n$-spherical theories.

We use without specifications  notions, notations and results on
arities of theories \cite{aafot, ananat}, on circular orders
\cite{KulMac, Kulp3, Kulp4} and on distributions of countable
models of complete theories \cite{CCMCT}.

\section{Preliminaries}\label{sec2}

Recall a series of notions related to arities and aritizabilities
of theories.

{\bf Definition} \cite{ZPS}. A theory $T$ is said to be
\emph{$\Delta$-based}\index{Theory!$\Delta$-based}, where $\Delta$
is some set of~formulae without parameters, if any formula of $T$
is equivalent in $T$ to a~Boolean combination of formulae
in~$\Delta$.

For $\Delta$-based theories $T$, it is also said that $T$ has {\em
quantifier elimination}\index{Elimination of quantifiers} or {\em
quantifier reduction}\index{Reduction of quantifiers} up to
$\Delta$.

\medskip
{\bf Definition} \cite{CCMCT, ZPS}. {\rm Let $\Delta$ be a set of
formulae of a theory $T$, and $p(\bar{x})$ a type of $T$ lying in
$S(T)$. The type $p(\bar{x})$ is said to be
\emph{$\Delta$-based}\index{Type!$\Delta$-based} if $p(\bar{x})$
is isolated by a set of formulas $\varphi^\delta\in p$, where
$\varphi\in\Delta$, $\delta\in\{0,1\}$.}

\medskip
The following lemma, being a corollary of Compactness Theorem,
noticed in \cite{ZPS}.

\begin{lemma}\label{lem01} A  theory  $T$  is
$\Delta$-based  if and only if,  for  any  tuple~$\bar{a}$ of any
{\rm (}some{\rm )} weakly saturated model of $T$, the type ${\rm
tp}(\bar{a})$ is $\Delta$-based.
\end{lemma}

Lemma \ref{lem01} allows to reduce, in general, possibilities of
$\Delta$-basedness of a theory $T$ to the $\Delta$-deducibility of
its types, and vice versa.

Below we consider special forms of $\Delta$-basedness.

\medskip
{\bf Definition} \cite{aafot}. An elementary theory $T$ is called
{\em unary}, or {\em $1$-ary}, if any $T$-formula
$\varphi(\overline{x})$ is $T$-equivalent to a Boolean combination
of $T$-formulas, each of which is of one free variable, and of
formulas of form $x\approx y$.

For a natural number $n\geq 1$, a formula $\varphi(\overline{x})$
of a theory $T$ is called {\em $n$-ary}, or an {\em $n$-formula},
if $\varphi(\overline{x})$ is $T$-equivalent to a Boolean
combination of $T$-formulas, each of which is of $n$ free
variables.

For a natural number $n\geq 2$, an elementary theory $T$ is called
{\em $n$-ary}, or an {\em $n$-theory}, if any $T$-formula
$\varphi(\overline{x})$ is $n$-ary.

A theory $T$ is called {\em binary} if $T$ is $2$-ary, it is
called {\em ternary} if $T$ is $3$-ary, etc.

We will admit the case $n=0$ for $n$-formulae
$\varphi(\overline{x})$. In such a case $\varphi(\overline{x})$ is
just $T$-equivalent to a sentence
$\forall\overline{x}\varphi(\overline{x})$.

If $T$ is a theory such that $T$ is $n$-ary and not $(n-1)$-ary
then the value $n$ is called the arity of $T$ and it is denoted by
${\rm ar}(T)$. If $T$ does not have any arity we put ${\rm
ar}(T)=\infty$.

Similarly, for a formula $\varphi$ of a theory $T$ we denote by
${\rm ar}_T(\varphi)$ the natural value $n$ if $\varphi$ is
$n$-ary and not $(n-1)$-ary. If a theory $T$ is fixed we write
${\rm ar}(\varphi)$ instead of ${\rm ar}_T(\varphi)$.

Clearly, ${\rm ar}(\varphi)\leq\vert{\rm FV}(\varphi)\vert$, where
${\rm FV}(\varphi)$ is the set of free variables of the formula
$\varphi$.

\section{$n$-spherical orders}\label{sec3}

Recall \cite{KulMac, Kulp3, Kulp4} that a {\em circular}, or {\em
cyclic}  order relation is
described by a ternary relation $K_3$ satisfying the following conditions:\\
 (co1) $\forall x\forall y \forall z (K_3(x,y,z)\to K_3(y,z,x));$\\
 (co2) $\forall x\forall y \forall z (K_3(x,y,z)\land K_3(y,x,z)
 \leftrightarrow x\approx y \lor y\approx z \lor z\approx x);$\\
 (co3) $\forall x\forall y \forall z(K_3(x,y,z)\to \forall t[K_3(x,y,t)
 \lor K_3(t,y,z)]);$\\
 (co4) $\forall x\forall y \forall z (K_3(x,y,z)\lor K_3(y,x,z)).$

Clearly, ${\rm ar}(K_3(x,y,z))=3$ if the relation has at least
three element domain. Hence, theories with infinite circular order
relations are at least $3$-ary.

The following generalization of circular order produces a {\em
$n$-ball}, or {\em $n$-spherical}, or {\em $n$-circular} order
relation \cite{aafot, ananat}, for $n\geq 4$, which is described
by a $n$-ary relation
$K_n$ satisfying the following conditions:\\
 (nso1) $\forall x_1,\ldots,x_n (K_n(x_1,x_2,\ldots,x_n)\to K_n(x_2,\ldots,x_n,x_1));$\\
 (nso2) $\forall x_1,\ldots,x_n
 \biggm((K_n(x_1,\ldots,x_i,\ldots,x_j,\ldots,x_n)\land$
 $$\land K_n(x_1,\ldots,x_j,\ldots,x_i,\ldots,x_n))\leftrightarrow\bigvee\limits_{1\leq k<l\leq n} x_k\approx x_l\biggm)$$ for any $1\leq i<j\leq n$;\\
 (nso3) $\forall x_1,\ldots,x_n\Biggm(K_n(x_1,\ldots,x_n)\to$
 $$\to\forall t\left(\bigvee\limits_{i=1}^nK_n(x_1,\ldots,x_{i-1},t,x_{i+1},\ldots,x_n)\right)\Biggm);$$\\
 (nso4) $\forall x_1,\ldots,x_n (K_n(x_1,\ldots,x_i,\ldots,x_j,\ldots,x_n)\lor$
 $$\lor K_n(x_1,\ldots,x_j,\ldots,x_i,\ldots,x_n)),\,\, 1\leq i<j\leq
 n.$$

Structures $\mathcal{A}=\langle A,K_n\rangle$ with $n$-spherical
orders will be called {\em $n$-spherical orders}, too.

\begin{note} {\rm The axioms above are valid for $n=2$ producing a linear
order $K_2$. The only case $n=2$ can admit endpoints.}
\end{note}

Now we consider a series of models illustrating $n$-spherical
orders.

\begin{example}\label{ex12} {\rm Recall that a {\em directed tetrahedron}, or {\em
dirtetrahedron} \cite{MP} is a $4$-tuple $(A_1,A_2,A_3,$ $A_4)$ of
points in $\mathbb R^3$ which do not belong to a common plane. For
the dirtetrahedron $(A_1,A_2,A_3,A_4)$ we denote by $K^4_4$ the
closure of the set of even permutations of $(A_1,A_2,A_3,A_4)$ by
cyclic permutations united with the set $I_4$ of $4$-tuples
$(A_i,A_j,A_k,A_l)$ having at least two equal points. Thus $K^4_4$
is generated by three dirtetrahedrons $$(A_1,A_2,A_3,A_4),
(A_1,A_3,A_4,A_2), (A_1,A_4,A_2,A_3)$$ using cyclic permutations
and $I_4$.

Clearly, $K^4_4$ satisfies the axioms (nso1)--(nso4) for $n=4$,
producing a $4$-spherical order. It contains all $4$-tuple which
do not produce odd permutations of $(A_1,A_2,A_3,A_4)$, i.e.,
$\vert K^4_4 \vert =4^4-\frac{4!}{2}=4(4^3-3)=244$.

\begin{figure}
\begin{center}
\setlength{\unitlength}{1.1mm}
\begin{picture}(45,45)(4.5,8)

\qbezier(45.0, 30.0)(45.0, 38.2843) (39.1421, 44.1421)
\qbezier(39.1421, 44.1421)(33.2843, 50.0) (25.0, 50.0)
\qbezier(25.0, 50.0)(16.7157, 50.0) (10.8579, 44.1421)
\qbezier(10.8579, 44.1421)(5.0, 38.2843) (5.0, 30.0) \qbezier(5.0,
30.0)(5.0, 21.7157) (10.8579, 15.8579) \qbezier(10.8579,
15.8579)(16.7157, 10.0) (25.0, 10.0) \qbezier(25.0, 10.0)(33.2843,
10.0) (39.1421, 15.8579) \qbezier(39.1421, 15.8579)(45.0, 21.7157)
(45.0, 30.0)

\put(45,30){\makebox(0,0)[cc]{$\bullet$}}
\put(5,30){\makebox(0,0)[cc]{$\bullet$}}
\put(25,50){\makebox(0,0)[cc]{$\bullet$}}
\put(25,10){\makebox(0,0)[cc]{$\bullet$}}

\put(48,30){\makebox(0,0)[cc]{$a_1$}}
\put(25,52.5){\makebox(0,0)[cc]{$a_2$}}
\put(2.5,30){\makebox(0,0)[cc]{$a_3$}}
\put(25.7,7.5){\makebox(0,0)[cc]{$a_4$}}

\put(25.6,49.99){\vector(-4,0){0.5}}
\put(5,31){\vector(0,-4){0.5}} \put(24.2,10.01){\vector(4,0){0.5}}
\put(45,29){\vector(0,4){0.5}}

\put(45,30){\vector(-4,0){39.5}} \put(5,30){\vector(1,-1){19.5}}
\put(25,10){\vector(0,1){39.5}} \put(25,50){\vector(1,-1){19.5}}

\qbezier(45,30)(35,28)(25,10) \qbezier(25,10)(17,30)(25,50)
\qbezier(25,50)(13,34)(5,30) \qbezier(5,30)(25,38)(45,30)

\put(25.5,11){\vector(-1,-1){0.5}}
\put(24.05,46.85){\vector(1,4){0.5}}
\put(6.85,31.2){\vector(-3,-2){0.5}}
\put(43.3,30.6){\vector(3,-1){0.5}}

\end{picture}
\end{center}
\hspace{54mm}\parbox[t]{0.3\textwidth}\caption{} \label{fig111}
\end{figure}

Replacing points $A_1,A_2,A_3,A_4$ by elements $a_1,a_2,a_3,a_4$,
respectively, we obtain an isomorphic copy for $K^4_4$ with
generating $4$-tuples represented in Figure \ref{fig111}. }
\end{example}

\begin{example}\label{ex13} {\rm
Similarly Example \ref{ex12}, taking a $n$-tuple $(A_1,A_2,\ldots,
A_n)$ of points in $\mathbb R^{n-1}$, $n\geq 5$, which do not
belong to a common hyperplane we obtain a {\em directed
$n$-hedron}. We denote by $K^n_n$ the closure of the set of even
permutations of $(A_1,A_2,\ldots, A_n)$ by cyclic permutations
united with the set $I_n$ of $n$-tuples $(A_{i_1},A_{i_2},\ldots,
A_{i_n})$ having at least two equal points.

Clearly, $K^n_n$ satisfies the axioms (nso1)--(nso4) in the
general case, producing a $n$-spherical order.

By the definition the relations $K^n_n$ are generated from the
tuple $(A_1,A_2,\ldots,A_n)$ by its even permutations, cyclic
permutations and identifications of coordinates.

Similarly to $K^4_4$ we obtain $\vert K^n_n \vert
=n^n-\frac{n!}{2}$. This formula is valid for $2$-element linear
orders and $3$-element circular orders. Thus we have $\vert
K^2_2\vert =3$, $\vert K^3_3\vert =24$, $\vert K^4_4\vert =244$,
$\vert K^5_5\vert =3065$, etc.

For instance, $K^5_5$ is generated by twelve directed pentahedrons
$$(A_1,A_2,A_3,A_4,A_5), (A_1,A_2,A_4,A_5,A_3), (A_1,A_2,A_5,A_3,A_4),$$
$$(A_1,A_3,A_2,A_5,A_4), (A_1,A_3,A_4,A_2,A_5), (A_1,A_3,A_5,A_4,A_2).$$
$$(A_1,A_4,A_2,A_3,A_5), (A_1,A_4,A_3,A_5,A_2), (A_1,A_4,A_5,A_2,A_3),$$
$$(A_1,A_5,A_2,A_4,A_3), (A_1,A_5,A_3,A_2,A_4), (A_1,A_5,A_4,A_3,A_2)$$
using cyclic permutations and $I_5$.}
\end{example}

\begin{note}\label{rem131} {\rm Constructing $n$-spherical orders
on $n$-element sets we start to use even permutations since $n=4$
because for $n=2$ there are no nonidentical even permutations and
for $n=3$ these permutations are reduced to cyclic permutations.}
\end{note}

We denote by $S_n:=S(A_1,A_2,\ldots,A_n)$ the unique
$(n-1)$-dimensional sphere containing the points
$A_1,A_2,\ldots,A_n$ in Example \ref{ex13}, and by
$H_n:=H(A_1,A_2,\ldots,A_n)$ the union of hyperplanes containing
$(n-1)$-element subsets of $\{A_1,A_2,\ldots,A_n\}$.

\begin{example}\label{ex14} {\rm
Any relation $K^n_n$ can be extended till a $n$-spherical order
$K^{n+k}_n$ by new points $B_1,\ldots,B_k\in S_n\setminus H_n$ and
generating $n$-tuples
$(A_{i_1},\ldots,A_{i_r},B_{j_1},\ldots,B_{j_s})$ with $i_1<\ldots
<i_r$, $j_1<\ldots <j_s$ with respect to even permutations, cyclic
permutations and identifications of coordinates.

The structures with the relations $K^{n+k}_n$ generalize the
structures with the relations $K^n_n$ in Example \ref{ex13}. Since
$K^{n+k}_n$ admits even permutations only and there are
$A_{n+k}^n$ $n$-tuples with distinct coordinates we obtain
$$\vert K^{n+k}_n\vert =(n+k)^n-\frac{A_{n+k}^n}{2}.$$

In particular, we have $\vert K^3_2 \vert =3^2-\frac{3!}{2}=6$,
$\vert K^4_2\vert =4^2-\frac{4!}{4}=10$, $\vert K^5_2 \vert
=5^2-\frac{5!}{12}=15$, etc.; $\vert K^4_3\vert
=4^3-\frac{4!}{2}=52$, $\vert K^5_3\vert =5^3-\frac{5!}{4}=95$,
$\vert K^6_3\vert =6^3-\frac{6!}{12}=156$, etc.; $\vert K^5_4\vert
=5^4-\frac{5!}{2}=565$, $\vert K^6_4\vert =6^4-\frac{6!}{4}=1116$,
$\vert K^7_4\vert =7^4-\frac{7!}{12}=1981$, etc. }
\end{example}

\begin{note}\label{rem141} {\rm The construction of finite $n$-spherical orders
shows that these orders of given finite cardinality are unique up
to isomorphism.}
\end{note}

\begin{note}\label{rem15}
{\rm The process of adding $B_1,\ldots,B_k$ in Example \ref{ex14}
can be continued arbitrarily many (unlimit and limit) steps
obtaining $\lambda$-element $n$-spherical orders $K^\lambda_n$ on
$S_n$, for $\lambda$ with $n\leq\lambda\leq 2^\omega$, and on
appropriate elementary extensions $S^\lambda_n$ of $S_n$, for
$\lambda> 2^\omega$. All these orders on spheres are called {\em
spherical models}.

The definition of $K^\lambda_n$ implies that this relation does
not reduced to Boolean combinations of relations with less
arities, then the arity of the formula
$K^\lambda_n(x_1,\ldots,x_n)$ equals $n$: ${\rm
ar}(K^\lambda_n(x_1,\ldots,x_n))=n$.

The formula witnessing the identification of coordinates for the
formula $K_n(x_1,\ldots,x_n)$ of $n$-spherical order will be
denoted by $I_n(x_1,\ldots,x_n)$. Clearly, ${\rm
ar}(I_n(x_1,\ldots,x_n))=1$, since $I_n(x_1,\ldots,x_n)$ is
equivalent to the formula $$\bigvee\limits_{1\leq i<j\leq n}
\left(x_i\approx
x_j\wedge\bigwedge\limits_{k\notin\{i,j\}}x_k\approx x_k\right).$$
Thus, ${\rm ar}(K_n(x_1,\ldots,x_n)\wedge\neg
I_n(x_1,\ldots,x_n))=n$.}
\end{note}

\section{Countably categorical $n$-spherical orders}\label{sec4}

{\bf Definition.} A $n$-spherical order $K_n$, $n\geq 2$, is
called {\em dense} if it contains at least two elements and for
each $(a_1,a_2,a_3,\ldots,a_n)\in K_n$ with $a_1\ne a_2$ there is
$b\notin\{a_1,a_2,\ldots,a_n\}$ with
$$
\models K_n(a_1,b,a_3,\ldots,a_n)\wedge K_n(b,a_2,a_3,\ldots,a_n).
$$

Clearly, dense $n$-spherical orders $K_n$ witness the strict order
property producing unstable structures $\langle A,K_n\rangle$.

The following theorem generalizes the known result on countable
categoricity of dense linear orders \cite[Proposition
3.1.7]{ErPa}.

\begin{theorem}\label{th_dense_isom}
If $\mathcal{A}$ and $\mathcal{B}$ are countable dense
$n$-spherical orders, $n\geq 2$, without endpoints for $n=2$, then
$\mathcal{A}\simeq\mathcal{B}$.
\end{theorem}

Proof. The case $n=2$ is considered in \cite[Proposition
3.1.7]{ErPa}. Therefore we assume that $n\geq 3$, and for this
case we slightly modify the arguments for that proposition.

Let $A=\{a_m\mid m\in\omega\}$, $B=\{b_m\mid m\in\omega\}$. We
consider the set $G$ consisting of maps $g\mbox{\rm : }
A_1\rightarrow B_1$ satisfying the following conditions:

1) $A_1$ and $B_1$ are finite subsets of $A$ and $B$,
respectively;

2) $g\mbox{\rm : }\mathcal{A}(A_1)\iso\mathcal{B}(B_1)$ if
$A_1\neq\emptyset$;

3) if $\vert A_1\vert =2m>0$ then $\{a_0,\dots,a_{m-1}\}\subseteq
A_1$ and $\{b_0,\dots,b_{m-1}\}\subseteq B_1$;

4) if $\vert A_1 \vert =2m+1$ then $\{a_0,\dots,a_m\}\subseteq
A_1$ and for $m>0$, $\{b_0,\dots,b_{m-1}\}\subseteq B_1$.

Since $\emptyset\in G$ then $G\neq\emptyset$. Let $g\mbox{\rm : }
A_1\rightarrow B_1$ belong to $G$ and $\vert A_1\vert =2n$. In
view of 3) there is $a\in A\setminus A_1$ with
$\{a_0,\dots,a_m\}\subseteq A_1\cup\{a\}$. Now we find an element
$b\in B\setminus B_1$ such that
$K_n(b,g(c_2),\ldots,g(c_n))\Leftrightarrow K_n(a,c_2,\ldots,c_n)$
for all $c_2,\ldots,c_n\in A_1$. The element $b$ exists since
$\mathcal{B}$ is dense and the condition 2) for $g$ holds.
Clearly, $g\cup\{\langle a,b \rangle\}\in G$. If $\vert A_1\vert
=2n+1$ then we replace $\mathcal{A}$ and $\mathcal{B}$ each other
and find a pair $\langle a,b \rangle\notin g$ with $g\cup\{\langle
a,b \rangle\}\in G$. Thus the partial order $\langle
G,\subseteq\rangle$ does not have maximal elements. Hence $G$
contains an infinite chain $X\subseteq G$. In view of 2)--4) the
union $\bigcup X$ is an isomorphism of $\mathcal{A}$ onto
$\mathcal{B}$. \endproof

\medskip
Theorem \ref{th_dense_isom} immediately implies:

\begin{corollary}\label{cor_dcc}
For any dense $n$-spherical order $\mathcal{A}$ its theory ${\rm
Th}(\mathcal{A})$ is countably categorical.
\end{corollary}

\begin{note}\label{rem_ccso}
{\rm Extensions of dense $n$-spherical orders by finitely many new
discrete elements preserve the countable categoricity. It means
that there are finitely many possibilities for
$(a_1,a_2,a_3,\ldots,a_n)\in K_n$ such that $a_1\ne a_2$ and there
are finitely many $b\notin\{a_1,a_2\}$ with
$$
\models K_n(a_1,b,a_3,\ldots,a_n)\wedge K_n(b,a_2,a_3,\ldots,a_n).
$$
Similarly linear orders having infinite discrete parts loose the
countable categoricity obtaining a definable successor function.}
\end{note}

Similarly theories of dense linear orders, in view of Lemma
\ref{lem01} theories $T_n$ of dense $n$-spherical orders $\langle
A, K_n\rangle$, $n\geq 2$, admit the quantifier elimination since
complete types are forced by collections of quantifier free
formulas. Besides the theories $T_n$ are finitely axiomatizable.
Using Corollary \ref{cor_dcc} and the arguments for
\cite[Proposition 8.3.1]{ErPa} we obtain the following its
generalization:

\begin{theorem}\label{th_dec}
For any natural $n\geq 2$ the theory $T_n$ of dense $n$-spherical
order is decidable.
\end{theorem}

\section{Ehrenfeucht theories based on dense $n$-sphe\-ri\-cal orders}\label{sec5}

In this section we modify classical Ehrenfeucht examples on linear
dense orders \cite{CCMCT, Va} to Ehrenfeucht theories based on
dense $n$-sphe\-ri\-cal orders and producing arbitrarily finitely
many countable models.

For a theory $T$, we denote by
$I(T,\lambda)$\index{$I(T,\lambda)$} the number of pairwise
non-isomorphic models of $T$ in~a~power~$\lambda$. The value
$I(T,\lambda)$ is called the {\em $\lambda$-spectrum} of $T$. If
$\lambda=\omega$ then the $\lambda$-spectrum of $T$ is called the
{\em countable spectrum} of $T$.

\medskip
{\bf Definition} \cite{Mi2}. A theory $T$ is \emph{Ehrenfeucht} if
$1<I(T,\omega)<\omega,$ i.e., the countable spectrum of $T$ is
finite and greater than $1$.

\medskip
Let $\mathcal{A}_n=\langle A_n, K_n\rangle$ be a countable
dense $n$-spherical order, $n\geq 2$. Let $T^m_n$ be the theory of
a structure $\mathcal{M}^m_n$, formed from the structure
$\mathcal{A}_n$ by adding pairwise distinct constants $c_k$,
$k\in\omega$, such that $\models K_n(c_{i_1},\ldots,c_{i_n})$ for
any $i_1<\ldots<i_n$, and unary predicates $P_0,\ldots,P_{m-3}$
which form a~partition of the set $A_n$, with
$$\models\forall x_1,\ldots,x_n\:\Biggm(\bigwedge\limits_{1\leq j<k\leq n}\neg(x_j\approx x_k)\wedge K_n(x_1,\ldots,x_n)\to$$
$$\to\exists t\:\Biggm(\bigwedge\limits_{j=1}^n\neg(t\approx
x_j)\wedge K_n(x_1,t,x_3,\ldots,x_n)\wedge
K_n(t,x_2,x_3,\ldots,x_n)\wedge P_i(t)\Biggm)\Biggm),$$
$i=0,\ldots,m-3.$

We denote by $p_\infty(x)$ the type which is forced by the set of
formulae $K_n(c_{i_1},c_{i_2},\ldots,c_{i_{n-1}},x)$,
$i_1<i_2<\ldots<i_{n-1}$. The type $p_\infty(x)$ has $m-2$
completions $p^i_\infty(x)$ which are forced by
$p_\infty(x)\cup\{P_i(x)\}$, $i=0,\ldots,m-3.$

The theory $T^m_n$ has exactly $m$ pairwise non-isomorphic
countable models:

(a) a prime model $\mathcal{M}_m$ which omit the type
$p_\infty(x)$; ($\lim\limits_{k\to\infty}c_k=\infty$);

(b) prime models $\mathcal{M}^i_m$ over realizations of types
$p^i_\infty(x)$, $i=0,\ldots,m-3$, with a limit element for the
constants $c_k$, realizing $p^i_\infty(x)$;

(c) a saturated model $\overline{\mathcal{M}}_m$, without limit
elements for the constants $c_k$, realizing $p^i_\infty(x)$.

Replacing constants $c_k$ by unary predicates $U_k=\{c_k\}$ we
obtain theories $\widehat{T}^m_n$ instead of $T^m_n$ such that
$\widehat{T}^m_n$ are unary expansions of $T_n={\rm
Th}(\mathcal{A}_n)$, i.e. expansions of $T_n$ by unary predicates.
Clearly, this transformation preserves the number of pairwise
non-isomorphic countable models:
$$
I(\widehat{T}^m_n,\omega)=m=I(T^m_n,\omega).
$$

Corollary \ref{cor_dcc} and the arguments above producing values
for $ I(\widehat{T}^m_n,\omega)$ imply the following:

\begin{theorem}\label{th_dso_spec}
For any $m\in\omega\setminus\{0,3\}$, $n\in\omega\setminus\{0,1\}$
there is a unary expansion $\widehat{T}^m_n$ of the theory $T_n$
of dense $n$-spherical order $\mathcal{A}_n$ such that $
I(\widehat{T}^m_n,\omega)=m$.
\end{theorem}

The theories $T^m_2$ in Theorem \ref{th_dso_spec} are classical
Ehrenfeucht examples, and the theories $T^m_3$ are their
adaptations for dense circular orders.

Figure \ref{fig222a} illustrates a model of the theory $T^3_4$.
Constants $c_k$, $k\in\omega$, are situated on a spiral with a top
limit vertex $\ast$.

\begin{figure}
\begin{center}
\setlength{\unitlength}{1.1mm}
\begin{picture}(45,45)(4.5,11)

\qbezier(45.0, 30.0)(45.0, 38.2843) (39.1421, 44.1421)
\qbezier(39.1421, 44.1421)(33.2843, 50.0) (25.0, 50.0)
\qbezier(25.0, 50.0)(16.7157, 50.0) (10.8579, 44.1421)
\qbezier(10.8579, 44.1421)(5.0, 38.2843) (5.0, 30.0) \qbezier(5.0,
30.0)(5.0, 21.7157) (10.8579, 15.8579) \qbezier(10.8579,
15.8579)(16.7157, 10.0) (25.0, 10.0) \qbezier(25.0, 10.0)(33.2843,
10.0) (39.1421, 15.8579) \qbezier(39.1421, 15.8579)(45.0, 21.7157)
(45.0, 30.0)

\put(25,50){\makebox(0,0)[cc]{$\ast$}}
\put(20,15.25){\makebox(0,0)[cc]{$\bullet$}}
\put(20.5,13.2){\makebox(0,0)[cc]{$c_0$}}
\put(40,24){\makebox(0,0)[cc]{$\bullet$}}
\put(41,22){\makebox(0,0)[cc]{$c_1$}}
\put(41,39){\makebox(0,0)[cc]{$\bullet$}}
\put(40.9,37){\makebox(0,0)[cc]{$c_2$}}

\put(18,42.6){\makebox(0,0)[cc]{$\bullet$}}
\put(18,45){\makebox(0,0)[cc]{$c_3$}}

\qbezier(5,30)(25,18)(45,30) 

\qbezier(5,30)(6.5,30.95)(8,31.45)
\qbezier(10,32.2)(11.5,32.85)(13,33.2)
\qbezier(15,33.8)(16.5,34.2)(18,34.5)
\qbezier(20,34.9)(23,35.4)(24,35.3)
\qbezier(26,35.3)(27,35.4)(30,34.9)
\qbezier(32,34.5)(33.5,34.2)(35,33.8)
\qbezier(37,33.2)(38.5,32.85)(40,32.2)
\qbezier(42,31.45)(43.5,30.95)(45,30)

\qbezier(25,10)(13,30)(25,50)

\qbezier(25,10)(25.95,11.5)(26.45,13)
\qbezier(27.2,15)(28.0,17.5)(28.15,18)
\qbezier(28.7,20)(29.1,21.5)(29.4,23)
\qbezier(29.8,25)(30.3,28)(30.2,29)
\qbezier(30.2,31)(30.3,34)(29.8,35)
\qbezier(29.4,37)(29.1,38.5)(28.7,40)
\qbezier(28.15,42)(27.8,43.5)(27.2,45)
\qbezier(26.45,47)(25.95,48.5)(25,50)

\qbezier(15,15)(43,16)(44.1,35.8)
\qbezier(44.1,35.8)(43,39)(38,40.3)
\qbezier(36,40.8)(33,41.7)(30,42)
\qbezier(28,42.3)(25,42.8)(22,42.7)
\qbezier(20,42.7)(17,42.8)(14,42.6)
\qbezier(12,42.51)(11.5,42.55)(9.2,42.29)
\qbezier(9.2,42.29)(36.9,37.55)(34.5,47.5)

\put(32,47.6){\makebox(0,0)[cc]{$\ldots$}}

\end{picture}
\end{center}
\hspace{54mm}\parbox[t]{0.3\textwidth}\caption{} \label{fig222a}
\end{figure}

\section{Constant expansions of dense $n$-sphe\-ri\-cal orders and their countable spectra}\label{sec6}

Now we consider possibilities for distributions of countable
sequences of constants expanding $n$-spherical theories $T_n={\rm
Th}(\mathcal{A}_n)$, where $\mathcal{A}_n=\langle A_n,K_n\rangle$
are dense $n$-spherical orders, $n\geq 2$. These distributions
produce distributions of countable models of these expansions
$T\supseteq T_n$ and possibilities for values $I(T,\omega)$.

The distributions for $n=2$ are described in \cite{CCMCT, Ma}.

Recall \cite{PilSt} that a linearly ordered structure
$\mathcal{M}$ is {\em $o$-minimal}\index{Structure!$o$-minimal} if
any definable (with parameters) subset of $M$ is a finite union of
singletons and open intervals $(a,b)$, where $a\in
M\cup\{-\infty\}$, $b\in M\cup\{+\infty\}$. A theory $T$ is {\em
$o$-minimal}\index{Theory!$o$-minimal} if each model of $T$ is
$o$-minimal.

As examples of Ehrenfeucht $o$-minimal theories, we mention the
theories $T^1\rightleftharpoons{\rm Th}(({\mathbb
Q};<,c_n)_{n\in\omega})$ and $T^2\rightleftharpoons{\rm
Th}(({\mathbb Q};<,c_n,c'_n)_{n\in\omega})$, where $<$ is an
ordinary strict order on the set ${\mathbb Q}$ of rationals,
constants $c_n$ form a strictly increasing sequence, and constants
$c'_n$ form a strictly decreasing sequence, $c_n<c'_n$,
$n\in\omega$.

The theory $T^1$ is an Ehrenfeucht's example \cite{Va} with
$I(T^1,\omega)=3$. It has two almost prime models and one limit
model:

{\small $\bullet$} a prime model with empty set of realizations of
type $p(x)$ isolated by the set $\{c_n<x\mid n\in\omega\}$ of
formulas;

{\small $\bullet$} a prime model over a realization of the type
$p(x)$, with the least realization of that type;

{\small $\bullet$} one limit model over the type $p(x)$, with the
set of realizations of $p(x)$ forming an open convex set.

The Hasse diagram for the Rudin--Keisler preorder $\leq_{\rm RK}$
and values of the function ${\rm IL}$ of distributions of numbers
of limit models for $\sim_{\rm RK}$-classes of $T^1$ is
represented in Figure~\ref{fig3}. This Hasse diagram equals the
Hasse diagram for the theory $T^3_4$.

\begin{figure}[b]
\begin{center}
\unitlength 4mm
\begin{picture}(5,6.5)(-2.5,1.0)
{\footnotesize\put(2,2.5){\line(0,5){5}}
\put(2,2.5){\makebox(0,0)[cc]{$\bullet$}}
\put(2,7.5){\makebox(0,0)[cc]{$\bullet$}}
\put(2,2.5){\circle{0.6}} \put(2,7.5){\circle{0.6}}
\put(3,2.5){\makebox(0,0)[cr]{$0$}}
\put(3,7.5){\makebox(0,0)[cr]{$1$}} }
\end{picture}
\hfill \unitlength 5mm
\begin{picture}(6,7)(4,1.7)
{\footnotesize \put(7.5,2.5){\line(0,5){5}}
\put(7.5,2.5){\makebox(0,0)[cc]{$\bullet$}}
\put(7.5,5){\makebox(0,0)[cc]{$\bullet$}}
\put(7.5,7.5){\makebox(0,0)[cc]{$\bullet$}}
\put(7.5,2.5){\circle{0.6}} \put(7.5,5){\circle{0.6}}
\put(7.5,7.5){\circle{0.6}} \put(8.5,2.5){\makebox(0,0)[cr]{$0$}}
\put(8.5,5){\makebox(0,0)[cr]{$0$}}
\put(8.5,7.5){\makebox(0,0)[cr]{$3$}} }
\end{picture}
\hfill \unitlength 5mm
\begin{picture}(6,7)(4,1.7)
{\footnotesize \put(6.0,2.5){\line(0,5){5}}
\put(6.0,2.5){\makebox(0,0)[cc]{$\bullet$}}
\put(6.0,5){\makebox(0,0)[cc]{$\bullet$}}
\put(6.0,7.5){\makebox(0,0)[cc]{$\bullet$}}
\put(6.0,2.5){\circle{0.6}} \put(6.0,5){\circle{0.6}}
\put(6.0,7.5){\circle{0.6}} \put(7.0,2.5){\makebox(0,0)[cr]{$0$}}
\put(7.0,5){\makebox(0,0)[cr]{$0$}}
\put(8.3,7.6){\makebox(0,0)[cr]{$2^k-1$}} }
\end{picture}

\end{center}
\hspace{-5mm}\parbox[t]{0.3\textwidth}{\caption{}\label{fig3}}
\hfill
\hspace{-66mm}\parbox[t]{0.47\textwidth}{\caption{}\label{fig4}}
\hfill
\hspace{-54mm}\parbox[t]{0.2\textwidth}{\caption{}\label{fig5}}

\end{figure}

The theory $T^2$ has six pairwise non-isomorphic countable models:

{\small $\bullet$} a prime model with empty set of realizations of
type $p(x)$ isolated by the set $\{c_n<x\mid
n\in\omega\}\cup\{x<c'_n\mid n\in\omega\}$;

{\small $\bullet$} a prime model over a realization of $p(x)$,
with a unique realization of this type;

{\small $\bullet$} a prime model over a realization of type
$q(x,y)$ isolated by the set $p(x)\cup p(y)\cup\{x<y\}$; here the
set of realizations of $p(x)$ forms a closed interval $[a,b]$;

{\small $\bullet$} three limit models over the type $q(x,y)$, in
which the sets of realizations of $q(x,y)$ are convex sets of
forms $(a,b]$, $[a,b)$, $(a,b)$ respectively.

In Figure \ref{fig4} we represent the Hasse diagram
of~Rudin--Keisler preorders $\leq_{\rm RK}$ and values of
distribution functions ${\rm IL}$ of numbers of limit models on
$\sim_{\rm RK}$-equivalence classes for the theory $T^2$.

The following theorem shows that the number of countable models of
Ehrenfeucht $o$-minimal theories is exhausted by combinations of
these numbers for the theories $T^1$ and $T^2$.

\begin{theorem}\label{th1Ma} {\rm \cite{Ma}}
Let $T$ be an $o$-minimal theory in a countable language. Then
either $T$ has $2^{\omega}$ countable models or $T$ has exactly
$3^r\cdot 6^s$ countable models, where $r$ and $s$ are natural
numbers. Moreover, for any $r,s\in\omega$ there is an $o$-minimal
theory $T$ with exactly $3^r\cdot 6^s$ countable models.
\end{theorem}

Theorem \ref{th1Ma} immediately implies:

\begin{corollary}\label{co2Ma}
Let $T$ be a countable constant expansion of the $2$-spherical
theory $T_2$. Then either $T$ has $2^{\omega}$ countable models or
$T$ has exactly $3^r\cdot 6^s$ countable models, where $r$ and $s$
are natural numbers. Moreover, for any $r,s\in\omega$ there is an
aforesaid theory $T$ with exactly $3^r\cdot 6^s$ countable models.
\end{corollary}

The arguments for the $2$-spherical theory $T_2$ are valid for the
$3$-spherical theory $T_3$ obtaining the following:

\begin{corollary}\label{co3Ma}
Let $T$ be a countable constant expansion of the $3$-spherical
theory $T_3$. Then either $T$ has $2^{\omega}$ countable models or
$T$ has exactly $3^r\cdot 6^s$ countable models, where $r$ and $s$
are natural numbers. Moreover, for any $r,s\in\omega$ there is an
aforesaid theory $T$ with exactly $3^r\cdot 6^s$ countable models.
\end{corollary}

Considering the dense $n$-spherical Ehrenfeucht theories $T_n$,
$n\geq 4$, we obtain both the possibilities $3^k\cdot 6^s$ of
countable models, with appropriate sequences of constants, and the
following new possibilities.

Considering a consistent nonisolated set $p(x)$ of formulae
$$K^\delta_n(c_{j_1},\ldots,c_{j_{i-1}},x,c_{j_{i+1}},\ldots,c_{j_n})\wedge\bigwedge\limits_{k\ne
i}\neg x\approx c_{j_k},$$
$\delta=\delta_{j_1,\ldots,j_{i-1},x,j_{i+1},\ldots,j_n}\in\{0,1\}$,
with pairwise distinct
$c_{j_1},\ldots,c_{j_{i-1}},c_{j_{i+1}},\ldots,c_{j_n}$, we
observe that each such an additional formula divides the domain of
previous ones into two parts with respect to a $(n-1)$-dimensional
plane containing points
$c_{j_1},\ldots,c_{j_{i-1}},c_{j_{i+1}},\ldots,c_{j_n}$. In view
of quantifier elimination for $T_n$ any consistent set $p(x)$
forces a complete type. Thus a limit part of domain being a set of
solutions for $p(x)$ is defined by the limits for
$c_{j_1},\ldots,c_{j_{i-1}},c_{j_{i+1}},\ldots,c_{j_n}$.

We have the following possibilities for these limits:

{\small $\bullet$} $p(x)$ is omitted;

{\small $\bullet$} $p(x)$ has a unique realization $a_\infty$
being a common limit of the sequences;

{\small $\bullet$} $p(x)$ has infinitely many realizations
independently including / not including $n-1$ limits with respect
to coordinates of
$(c_{j_1},\ldots,c_{j_{i-1}},c_{j_{i+1}},\ldots,c_{j_n})$.

Hence there are $2^{n-1}+2$ possibilities for countable models,
where 3 models are prime over finite sets (prime over $\emptyset$,
prime over $\{a_\infty\}$, and prime over limits for
$c_{j_1},\ldots,c_{j_{i-1}},c_{j_{i+1}},\ldots,c_{j_n}$) and
$2^{n-1}-1$ limit models.

If some constants in
$c_{j_1},\ldots,c_{j_{i-1}},c_{j_{i+1}},\ldots,c_{j_n}$  are
fixed, with $k>1$ independent moving sequences the total number of
countable models related to $p(x)$ equals $2^k+2$ including $3$
almost prime models and $2^k-1$ limit models, see Figure
\ref{fig5}.

Since distinct types $p(x)$ are independent the total number of
possibilities for countable models is obtained by the
multiplications of values $2^k+2$ for various $k\in
n\setminus\{1\}$, if there are finitely many nonisolated
$1$-types, and there are $2^{\omega}$ countable models otherwise.
Thus we have the following:

\begin{theorem}\label{thsphEr}
Let $T$ be a countable constant expansion of the $n$-spherical
theory $T_n$, $n\geq 3$. Then either $T$ has $2^{\omega}$
countable models or $T$ has exactly $\prod\limits_{k\in
n\setminus\{1\}}(2^k+2)^{r_k}$ countable models, where $r_k$ are
natural numbers. Moreover, for any $r_0,\ldots,r_{n-1}\in\omega$
there is an aforesaid theory $T$ with exactly $\prod\limits_{k\in
n\setminus\{1\}}(2^k+2)^{r_k}$ countable models.
\end{theorem}

Theorem \ref{thsphEr} confirms the Vaught conjecture for countable
constant expansions $T$ of $n$-spherical theories $T_n$. In
particular, $$\mbox{ either }I(T,\omega)=2^\omega\mbox{ or }
I(T,\omega)=3^{r_1}\cdot 6^{r_2}\cdot 10^{r_3}$$ for $T\supseteq
T_4$, $$\mbox{ either }I(T,\omega)=2^\omega\mbox{ or }
I(T,\omega)=3^{r_1}\cdot 6^{r_2}\cdot 10^{r_3}\cdot 18^{r_4}$$ for
$T\supseteq T_5$, etc.

Similarly to the quite $o$-minimal Ehrenfeucht theories
distributions of countable models of constant expansions of
$n$-spherical theories are given by Hasse diagrams for disjoint
unions of theories \cite{KulSudRK, RKDU}, based on the diagrams
represented in Figures \ref{fig3}--\ref{fig5}.

\section{Conclusion}\label{sec7}

We studied spherical orders, which generalize known linear and
circular orders. Semantic and syntactic properties of spherical
orders and their elementary theories, including finite and dense
orders and their theories are investigated. It is shown that
theories of dense $n$-spherical orders are countably categorical
and decidable. The values for spectra of countable models of unary
expansions of $n$-spherical theories are described. The Vaught
conjecture is confirmed for countable constant expansions of dense
$n$-spherical theories. It would be natural to study various
modifications and kinds of spherical orders and their theories.

Beibut Sh. Kulpeshov

Kazakh-British Technical University, Almaty, Kazakhstan;

Novosibirsk State Technical University, Novosibirsk, Russia

E-mail: kulpesh@mail.ru \vskip 2mm

Sergey V. Sudoplatov

Sobolev Institute of Mathematics,

Novosibirsk State Technical
University, Novosibirsk, Russia 

E-mail: sudoplat@math.nsc.ru

\end{document}